\documentclass[10pt,reqno]{amsart}

\usepackage{amssymb,amsmath,epsfig}
\usepackage{amsfonts}
\usepackage{pst-plot}

\newtheorem*{theorem*}{Theorem}

\numberwithin{equation}{section}

\title{On identities involving the sixth order mock theta functions}

\date{\today}
\author{Jeremy Lovejoy}
\address{CNRS, LIAFA, Universit\'e Denis Diderot - Paris 7,
Case 7014, 75205 Paris Cedex 13, FRANCE}
\email{lovejoy@liafa.jussieu.fr}
\subjclass[2000]{33D15}

\begin{document}
\begin{abstract}
We present $q$-series proofs of four identities involving sixth order mock theta functions from Ramanujan's lost notebook.  We also show how Ramanujan's identities can be used to give a quick proof of four sixth order identities of Berndt and Chan.  
\end{abstract}

\maketitle

\section{Ramanujan's sixth order identities}
The last four identities on p. 13 of Ramanujan's lost notebook \cite{Ra1} may be written as 
\begin{eqnarray}
q^{-1}\psi(q^2) + \rho(q) &=& (-q;q^2)_{\infty}^2(-q,-q^5,q^6;q^6)_{\infty}, \label{eq1} \\
\phi(q^2) + 2\sigma(q) &=& (-q;q^2)_{\infty}^2(-q^3,-q^3,q^6;q^6)_{\infty}, \label{eq2} \\
2\phi(q^2) - 2\mu(-q) &=& (-q;q^2)_{\infty}^2(-q^3,-q^3,q^6;q^6)_{\infty}, \label{eq3} \\ 
2q^{-1}\psi(q^2) + \lambda(-q) &=& (-q;q^2)_{\infty}^2(-q,-q^5,q^6;q^6)_{\infty}, \label{eq4}
\end{eqnarray}
where the ``sixth order" mock theta functions $\phi,\psi,\rho,\sigma,\lambda$, and $\mu$ are defined by 
\begin{eqnarray*}
\phi(q) &:=& \sum_{n \geq 0} \frac{(-1)^nq^{n^2}(q;q^2)_n}{(-q)_{2n}},  \\
\psi(q) &:=& \sum_{n \geq 0} \frac{(-1)^nq^{(n+1)^2}(q;q^2)_n}{(-q)_{2n+1}}, \\
\rho(q) &:=&  \sum_{n \geq 0} \frac{q^{n(n+1)/2}(-q)_n}{(q;q^2)_{n+1}}, \\
\sigma(q) &:=& \sum_{n \geq 0} \frac{q^{(n+1)(n+2)/2}(-q)_n}{(q;q^2)_{n+1}}, \\
\lambda(q) &:=& \sum_{n \geq 0} \frac{(-q)^n(q;q^2)_n}{(-q)_n}, \\
\mu(q) &:=& \sum_{n \geq 0} \frac{(-1)^n(q;q^2)_n}{(-q)_n}. 
\end{eqnarray*}

Here we use the usual $q$-series notation
\begin{equation*}
(a_1,a_2,\dots,a_k;q)_n := \prod_{j=0}^{n-1} (1-a_1q^j)(1-a_2q^j)\cdots(1-a_kq^j),
\end{equation*}
following the custom of dropping the ``$;q$" unless the base is something other than $q$.  The convergence of $\mu(q)$ (as well as \eqref{rr1} below) is in the ``Ces\`aro sense", which in this case means that the sum is obtained by averaging the limits of the even partial sums and the odd partial sums.

The identities \eqref{eq1} - \eqref{eq4} were proven by Andrews and Hickerson \cite{An-Hi1} by combining 
the Bailey pair method and the constant term method.  The first point of this note is that these identities 
follow immediately upon combining the $q$-series transformations
\begin{equation} \label{trans1}
\begin{aligned}
\frac{(-aq)_{\infty}}{(-q)_{\infty}}\sum_{n \geq 0} \frac{(x;q^2)_n(aq)_n(-q/x)^n}{(q^2;q^2)_n} &=
\frac{(-a^2q/x;q^2)_{\infty}}{(-q/x;q^2)_{\infty}}\sum_{n \geq 0} \frac{(a^2q^2;q^2)_{2n}(-1)^nq^{2n^2}}{(q^4;q^4)_n(-a^2q/x;q^2)_{2n+1}} \\ &- \frac{a(-a^2q^2/x;q^2)_{\infty}}{(-q^2/x;q^2)_{\infty}}\sum_{n \geq 1}\frac{(a^2q^2;q^2)_{n-1}(-q)^{n(n+1)/2}}{(-q;-q)_{n-1}(-a^2q^2/x;q^2)_n}
\end{aligned}
\end{equation} 
and
\begin{equation} \label{trans2}
\begin{aligned}
\frac{(-aq)_{\infty}}{(-q)_{\infty}}\sum_{n \geq 0} \frac{(x;q^2)_n(aq)_n(-q^2/x)^n}{(q^2;q^2)_n} &=
\frac{(-a^2q^2/x;q^2)_{\infty}}{(-q^2/x;q^2)_{\infty}}\sum_{n \geq 0}\frac{(a^2q^2;q^2)_{n}(-q)^{n(n+1)/2}}{(-q;-q)_{n}(-a^2q^2/x;q^2)_{n+1}}  \\ &+ \frac{a(-a^2q^3/x;q^2)_{\infty}}{(-q^3/x;q^2)_{\infty}}\sum_{n \geq 0}\frac{(a^2q^2;q^2)_{2n}(-1)^nq^{2n^2+4n+1}}{(q^4;q^4)_n(-a^2q^3/x;q^2)_{2n+1}}
\end{aligned}
\end{equation}
with the Rogers-Ramanujan type identities 
\begin{equation} \label{rr1}
2\sum_{n \geq 0} \frac{(q;q^2)_n(-1)^n}{(q)_n} = \frac{(q;q^2)_{\infty}(-q,-q^2,q^3;q^3)_{\infty}}{(q^2;q^2)_{\infty}}
\end{equation}
and
\begin{equation} \label{rr2}
\sum_{n \geq 0} \frac{(q;q^2)_n(-q)^n}{(q)_n} = \frac{(q;q^2)_{\infty}(-q^3,-q^3,q^3;q^3)_{\infty}}{(q^2;q^2)_{\infty}}.
\end{equation}
Indeed, taking $x=q$ and $a=\pm 1$ in \eqref{trans1} and \eqref{trans2} and appealing to \eqref{rr1} and \eqref{rr2} we obtain
\begin{eqnarray*}
\mu(q) &=& \frac{1}{2}\phi(q^2) - \sigma(-q), \\
(q;q^2)_{\infty}^3(-q,-q^2,q^3;q^3)_{\infty} &=& \phi(q^2) + 2\sigma(-q), \\
\lambda(q) &=& \rho(-q) + q^{-1}\psi(q^2), \\
(q;q^2)_{\infty}^3(-q^3,-q^3,q^3;q^3)_{\infty} &=& \rho(-q) - q^{-1}\psi(q^2), 
\end{eqnarray*}
and this easily yields \eqref{eq1} - \eqref{eq4}.  

The identities \eqref{rr1} and \eqref{rr2} come from letting $(a,b,c,d,e) \to (1,\infty,\infty,-\sqrt{q},\sqrt{q})$ and 
$(q,\infty,\infty,-\sqrt{q},\sqrt{q})$ in a limiting case of the Watson-Whipple transformation \cite[p. $62$, Ex. $2.22$]{Ga-Ra1},
\begin{equation*} 
\sum_{n \geq 0} \frac{(1-aq^{2n})(a,b,c,d,e)_n(-1)^nq^{n(n-1)/2}(aq)^{2n}}{(1-a)(q,aq/b,aq/c,aq/d,aq/e)_n(bcde)^n} = \frac{(aq,aq/de)_{\infty}}{(aq/d,aq/e)_{\infty}}\sum_{n \geq 0} \frac{(aq/bc,d,e)_n(aq/de)^n}{(q,aq/b,aq/c)_n}, 
\end{equation*}
and then applying the triple product identity \cite[p. $357$, Eq. (II.$28)$]{Ga-Ra1},
\begin{equation*} 
\sum_{n \in \mathbb{Z}} z^nq^{n^2} = (-zq,-q/z,q^2;q^2)_{\infty}. 
\end{equation*} 

As for \eqref{trans1} and \eqref{trans2}, when $x \to \infty$ these are two entries from Ramanujan's lost notebook, proven by Andrews and Berndt \cite[Entries $(1.4.6)$ and $(1.4.7)$]{An-Be1}.  Their proof generalizes in a straightforward way to give \eqref{trans1} and \eqref{trans2} as follows. 

\begin{proof}[Proof of \eqref{trans1} and \eqref{trans2}]
Two ``Heine-type" transformations \cite[Theorem $\text{A}_3$ and Theorem $\text{A}_1$]{An1} (or see \cite[Theorem $1.21$ and Theorem $1.22$]{An-Be1}) are required,
\begin{equation} \label{Heine1}
\sum_{n \geq 0} \frac{(a;q^2)_n(b)_{2n}t^n}{(q^2;q^2)_n(c)_{2n}} = \frac{(b)_{\infty}(at;q^2)_{\infty}}{(c)_{\infty}(t;q^2)_{\infty}}\sum_{n \geq 0}\frac{(c/b)_n(t;q^2)_nb^n}{(q)_n(at;q^2)_n}
\end{equation}
and
\begin{equation} \label{Heine2}
\begin{aligned}
\sum_{n \geq 0} \frac{(a;q^2)_n(b)_nt^n}{(q^2;q^2)_n(c)_n} &= \frac{(b)_{\infty}(at;q^2)_{\infty}}{(c)_{\infty}(t;q^2)_{\infty}}\sum_{n \geq 0} \frac{(c/b)_{2n}(t;q^2)_nb^{2n}}{(q)_{2n}(at;q^2)_n} \\
&+ \frac{(b)_{\infty}(atq;q^2)_{\infty}}{(c)_{\infty}(tq;q^2)_{\infty}}\sum_{n \geq 0} \frac{(c/b)_{2n+1}(tq;q^2)_nb^{2n+1}}{(q)_{2n+1}(atq;q^2)_n}.
\end{aligned}
\end{equation}

We begin with \eqref{trans1}.  First, setting $(q,a,b,c) = (q^2,q^2/t,a^2q^2,-a^2q^3/x)$ and then letting $t \to 0$ in \eqref{Heine1} we obtain
\begin{equation} \label{first}
\sum_{n \geq 0} \frac{(a^2q^2;q^2)_{2n}(-1)^nq^{2n^2}}{(q^4;q^4)_n(-a^2q/x;q^2)_{2n+1}} = \frac{(a^2q^2;q^2)_{\infty}(q^2;q^4)_{\infty}}{(-a^2q/x;q^2)_{\infty}}\sum_{n \geq 0} \frac{(-q/x;q^2)_na^{2n}q^{2n}}{(q^2;q^2)_{n}(q^2;q^4)_n}.
\end{equation}
Next letting $(a,b,c,t) \to (-q^2/x,0,-q^2,a^2q^2)$ and then setting $q=-q$ in \eqref{Heine1} gives
\begin{equation} \label{second}
\sum_{n \geq 0} \frac{(-q^2/x;q^2)_na^{2n+1}q^{2n+1}}{(q)_{2n+1}(-q^2;q^2)_n} = \frac{-a(-a^2q^2/x,-q;q^2)_{\infty}}{(a^2q^2;q^2)_{\infty}}\sum_{n \geq 1} \frac{(a^2q^2;q^2)_{n-1}(-q)^{n(n+1)/2}}{(-q;-q)_{n-1}(-a^2q^2/x;q^2)_{n}}.
\end{equation}
Then letting $(a,b,c,t) \to (x,aq,0,-q/x)$ in \eqref{Heine2} and multiplying both sides by $(-aq)_{\infty}/(-q)_{\infty}$
we have
\begin{equation} \label{third}
\begin{aligned}
\frac{(-aq)_{\infty}}{(-q)_{\infty}}\sum_{n \geq 0} \frac{(x;q^2)_n(aq)_n(-q/x)^n}{(q^2;q^2)_n} &= \frac{(a^2q^2;q^2)_{\infty}}{(-q^2,-q/x;q^2)_{\infty}}\sum_{n \geq 0} \frac{(-q/x;q^2)_na^{2n}q^{2n}}{(q)_{2n}(-q;q^2)_n} \\ &+\frac{(a^2q^2;q^2)_{\infty}}{(-q,-q^2/x;q^2)_{\infty}}\sum_{n \geq 0} \frac{(-q^2/x;q^2)_na^{2n+1}q^{2n+1}}{(q)_{2n+1}(-q^2;q^2)_n}.  
\end{aligned}
\end{equation}
Finally, applying \eqref{first} and \eqref{second} to the first and second terms on the right-hand side of \eqref{third} gives \eqref{trans1}.
 
The proof of \eqref{trans2} is similar so we just sketch it.   We let $(a,b,c,t) \to (x,aq,0,-q^2/x)$ in \eqref{Heine2}, multiply both sides by $(-aq)_{\infty}/(-q)_{\infty}$, then transform the first term on the right-hand side using the result of letting $(a,b,c,t) \to (-q^2/x,0,-q,a^2q^2)$ and then setting $q=-q$ in \eqref{Heine1} and transform the second term on the right hand side using the result of setting $(q,b,c,t) = (q^2,a^2q^2,-a^2q^5/x,q^6/a)$ and letting $a \to 0$ in \eqref{Heine1}.

\end{proof}

\section{Berndt and Chan's sixth order identities}
Recently Berndt and Chan \cite{Be-Ch1} defined two more sixth order mock theta functions, 
\begin{eqnarray*}
\phi_{-}(q) := \sum_{n \geq 1} \frac{(-q)_{2n-1}q^n}{(q;q^2)_n}, \\
\psi_{-}(q) := \sum_{n \geq 1} \frac{(-q)_{2n-2}q^n}{(q;q^2)_n}. 
\end{eqnarray*}
Using the same methods as Andrews and Hickerson, they proved the four identities
\begin{eqnarray}
-2q^{-1}\psi_{-}(q^2) + \rho(q) &=& (-q^2;q^2)_{\infty}^3(q^6,q^6,q^{12};q^{12})_{\infty}, \label{eq1bis} \\
-\phi_{-}(q^2) + \sigma(q) &=& q(-q^2;q^2)_{\infty}^2(-q^6,-q^6,q^6;q^6)_{\infty}, \label{eq2bis} \\
4\phi_{-}(q^2)  + 2\mu(q) &=& (-q;q^2)_{\infty}^2(-q^3,-q^3,q^6;q^6)_{\infty}, \label{eq3bis} \\ 
4q^{-1}\psi_{-}(q^2) + \lambda(q) &=& (-q;q^2)_{\infty}^3(q^3,q^9,q^{12};q^{12})_{\infty}. \label{eq4bis}
\end{eqnarray}
The second point of this note is that the identities of Berndt and Chan follow readily from those of Ramanujan.  Indeed, using equations \eqref{eq1} - \eqref{eq4} to eliminate $\rho,\sigma,\mu$, and $\lambda$ from equations \eqref{eq1bis} - \eqref{eq4bis}, we have the equivalent identities
\begin{eqnarray}
q^{-1}\psi(q^2) + 2q^{-1}\psi_{-}(q^2) &=&  (-q;q^2)_{\infty}^2(-q,-q^5,q^6;q^6)_{\infty} - (-q^2;q^2)_{\infty}^3(q^6,q^6,q^{12};q^{12})_{\infty}, \label{eq1ter} \\
\phi(q^2) + 2\phi_{-}(q^2) &=& (-q;q^2)_{\infty}^2(-q^3,-q^3,q^6;q^6)_{\infty} - 2q (-q^2;q^2)_{\infty}^2(-q^6,-q^6,q^6;q^6)_{\infty} \\
2\phi(q^2) + 4\phi_{-}(q^2) &=& (-q;q^2)_{\infty}^2(-q^3,-q^3,q^6;q^6)_{\infty} + (q;q^2)_{\infty}^2(q^3,q^3,q^6;q^6)_{\infty}, \\
2q^{-1}\psi(q^2) + 4q^{-1}\psi_{-}(q^2) &=& (-q;q^2)_{\infty}^2(-q,-q^5,q^6;q^6)_{\infty} - (q;q^2)_{\infty}^3(-q^3,-q^9,q^{12};q^{12})_{\infty}. \label{eq4ter}
\end{eqnarray}

The modularity of the left hand side in the case of $\phi + 2\phi_{-}$ follows from
an identity of Ramanujan \cite[Eq. $(3.4.2)$]{An-Be1},
\begin{equation*}
\sum_{n \geq 0} \frac{(-1)^n(q;q^2)_nq^{n^2}}{(-q)_{2n}} + 2\sum_{n \geq 1} \frac{(-q)_{2n-1}q^n}{(q;q^2)_n} = \frac{1}{(q)_{\infty}}\left(1+6\sum_{n \geq 0}\left(\frac{q^{6n+2}}{1-q^{6n+2}} - \frac{q^{6n+4}}{1-q^{6n+4}}\right)\right),
\end{equation*} 
combined with an identity of Lorenz \cite{Lo1} (or see \cite{Bo-Bo-Ga1}),
\begin{equation*}
\sum_{m,n \in \mathbb{Z}} q^{n^2+mn+m^2} = 1 + 6\sum_{n \geq 0}\left(\frac{q^{3n+1}}{1-q^{3n+1}} - \frac{q^{3n+2}}{1-q^{3n+2}}\right),
\end{equation*} 
while the modularity of $\psi+2\psi_{-}$ follows from the case $a \to -1$ of another identity of Ramanujan \cite[Eq. $(3.4.6)$]{An-Be1}, 
\begin{equation} \label{ramanother} 
\sum_{n \geq 0} \frac{(-aq;q^2)_na^{n+1}q^{(n+1)^2}}{(aq)_{2n+1}} - \sum_{n \geq 1}\frac{(1/a)_{2n-1}q^n}{(-q/a;q^2)_n} 
= \frac{q(q^2;q^2)_{\infty}(-a^3,-q^6/a^3,q^6;q^6)_{\infty}}{a(aq)_{\infty}(q;q^2)_{\infty}(-a,-q^2/a,-q/a,q^2;q^2)_{\infty}}.
\end{equation}
Now \eqref{eq1ter} - \eqref{eq4ter} follow from standard computational techniques for modular forms.  For example, using \eqref{ramanother} in \eqref{eq1ter}, setting $q=-q$, and doing a little rearranging we have that \eqref{eq1ter} is equivalent to
\begin{equation*}
\frac{-3\eta^4(24z)}{\eta^2(12z)} = \frac{\eta(8z)\eta^3(2z)\eta(24z)}{\eta(6z)\eta^2(4z)} - \frac{\eta^4(8z)}{\eta^2(4z)},
\end{equation*}
where $\eta(z) := q^{1/24}(q)_{\infty}$ and $q := e^{2 \pi i z}$.  This is an equality between holomorphic modular forms of weight $1$ on $\Gamma_0(24)$ (with a certain character), so its truth is established by verifying that the $q$-expansions of both sides agree up to $q^4$.  (Those unfamiliar with this method might consult \cite{On1}.)

\end{document}